\documentclass[12pt, reqno]{amsart}
\usepackage{amsfonts}
\usepackage{amssymb}
\usepackage{amsmath}
\usepackage{amsthm}
\usepackage{latexsym}
\usepackage{amsmath,amsfonts,amscd}

\let\qua\stdspace

\newtheorem{theorem}{Theorem}[section]
\newtheorem{lemma}[theorem]{Lemma}
\newtheorem{cor}[theorem]{Corollary}

\newcommand{\Tors}{\mathrm{Tors}}
\newcommand{\HF}{HF}

\newcommand{\Q}{\mathbb{Q}}
\newcommand{\R}{\mathbb{R}}

\newcommand{\Z}{\mathbb{Z}}

\newcommand{\OneHalf}{\frac{1}{2}}

\newcommand{\OneQuarter}{\frac{1}{4}}

\newcommand{\Zmod}[1]{\Z/{#1}\Z}

\newcommand{\Image}{\mathrm{Im}}

\newcommand{\cm}{\cdot}

\newcommand{\CDisk}{D}

\newcommand{\ModSWfour}{\mathcal{M}}
\newcommand{\ModFlow}{\ModSWfour}

\newcommand{\SpinC}{{\mathrm{Spin}}^c}

\newcommand\sgn{\mathrm{sgn}}

\newcommand\abuts\Rightarrow
\newcommand\Sym{\mathrm{Sym}}

\newcommand\Eul{\widehat{\chi}}

\newcommand\HFpRed{\HFp_{\red}}

\newcommand\chiTrunc{\chi^{\mathrm{trunc}}}

\newcommand{\ev}{\mathrm{ev}}

\newcommand\ModSphere{\ModFlow\left({\mathbb S}\longrightarrow
\Sym^{g-1}(\Sigma_{1})\times \Sym^2(\Sigma_{2})\right)}
\newcommand\ModSpheres\ModSphere

\newcommand\HFpred{\HFp_{\rm red}}

\newcommand{\red}{\mathrm{red}}

\newcommand\HFp{\HFb}

\newcommand\HFinf{HF^\infty}

\newcommand\HFb{HF^+}

\newcommand\UnparModSp{\widehat \ModSp}
\newcommand\UnparModFlow\UnparModSp
\newcommand\Mod\ModSp

\newcommand\PD{\mathrm{PD}}

\newcommand{\spinc}{\mathfrak s}

\newcommand{\spinct}{\mathfrak t}
\newcommand{\spincset}{\mathfrak T}
\newcommand{\spinca}{\mathfrak a}
\newcommand{\spincb}{\mathfrak b}
\newcommand{\spincc}{\mathfrak c}
\newcommand{\spincd}{\mathfrak d}
\newcommand\tSpinC{\underline{\SpinC}}
\newcommand\Yzero{Y_{p_0/q_0}}
\newcommand\Yone{Y_{p/q}}
\newcommand\Ytwo{Y_{p_2/q_2}}

\newcommand\ModMaps{\mathcal M}
\newcommand\ModSp\ModMaps

\newcommand\Dual{\mathcal D}
\newcommand\Duality\Dual

\begin{document}
\title[Surgery Formula for the Renormalized Euler Characteristic]{Surgery formula for the renormalized Euler characteristic of Heegaard Floer homology}

\author{Raif Rustamov}
\address{The Program in Applied and Computational Mathematics, Princeton University\\New Jersey 08540, USA}
\email{rustamov@princeton.edu} \keywords{Floer homology, surgery}
\begin{abstract}
We prove a surgery formula for renormalized Euler characteristic
of Ozsv\'ath and Szab\'o. The equality $\Eul=SW$ between this
Euler cahracteristic and the Seiberg-Witten invariant follows for
rational homology three-spheres.
\end{abstract}

\maketitle
\section {Introduction}

In \cite{HolDisk} and \cite{HolDiskFour} topological invariants
for closed oriented three manifolds and cobordisms between them
were defined by using a construction from symplectic geometry. The
resulting Floer homology package has many properties of a
topological quantum field theory.

Another such Floer homology package comes from Seiberg-Witten
theory \cite{KM}, \cite{MW}. Similarity of properties of the
Ozsv\'ath-Szab\'o and Seiberg-Witten theories and also
calculations heavily support the conjecture  that these invariants
are equivalent.

In this paper we will concentrate on a numerical invariant of
rational homology spheres obtained from the Heegaard Floer
homology package - the renormalized Euler characteristic, $\Eul$.
It is already known that for integral homology spheres $\Eul$ is
equal to Casson's invariant \cite{AbsGraded}, which is also the
case for the Seiberg-Witten invariant of integral homology spheres
\cite{Lim}. Calculations of \cite{Nemethi} push this equivalence
further to the Lens spaces and Seifert manifolds. Thus, it is
tempting to establish this equiality in its whole generality. To
this end we prove a surgery formula for $\Eul$. This formula and
several other properties of $\Eul$ and the related invariant
$\chiTrunc$ together fit into the framework of \cite{Nicola} to
give equivalence between $\Eul$ and the Reidemester-Turaev torsion
normalized by the Casson-Walker invariant. This also implies the
equality  $\Eul=SW$.

The organization of the paper is as follows: the required
preliminaries are presented in Section 2. The surgery theorem is
formulated and its applications are given in Section 3. The paper
finishes with the proof of the surgery formula in Section 4.

\medskip
\noindent{\bf{Acknowledgments}}\qua I am pleased to thank my
advisor Zolt\'an Szab\'o whose kind attitude and outstanding
teaching made this work very enjoyable. I would like to express my
gratitude to Paul Seymour for his support and encouragement.

\section{Preliminaries}

\noindent{\bf{Correction terms and Euler characteristics}}\qua Let
$Y$ be a rational homology sphere, $\spinct$ be a $\SpinC$
structure on it. We can consider Heegaard Floer homology group
$\HFp(Y,\spinct)$. This is a $\Q$ graded module over $Z[U]$. We
can also consider a simpler version, $\HFinf(Y,\spinct)$, for
which one can prove
\begin{equation}
\label{eq:HFinf} \HFinf(Y,\spinct) \cong Z[U,U^{-1}]
\end{equation}
 for each $\spinct$ structure. There is a natural $Z[U]$ equivariant map  $$\pi \colon\HFinf(Y,\spinct)\longrightarrow \HFp(Y,\spinct)$$ which is zero in sufficiently negative degrees and an isomorphism in all sufficiently positive degrees.
$\HFpred(Y, \spinct)$ is defined as
$$\HFpRed(Y, \spinct)=\HFp(Y,\spinct)/ \Image \pi.$$
Let $d(Y,\spinct)$ be the \emph{correction term} defined as the
minimal degree of any non-torsion class of $\HFp(Y, \spinct)$
lying in the image of $\pi$. Main object of our study, the
\emph{renormalized Euler characteristic} $\Eul(Y,\spinct)$ is
defined by
$$\Eul(Y, \spinct) = \chi(\HFpred(Y,\spinct)) - \frac{1}{2} d(Y,\spinct).$$

When $Y$ is a rational homology $S^1 \times S^2$ there is a
related numerical invariant $\chiTrunc$ as follows. Define
$\chiTrunc(Y, \spinct) = \chi(\HFp(Y, \spinct))$ for non-torsion
$\spinct$. If $\spinct$ is torsion then let $d(Y,\spinct)$ be the
minimal degree of any non-torsion class of $\HFp(Y, \spinct)$
coming from $\HFinf_{\ev}(Y,\spinct)$. The structure of $\HFinf$
for homology $S^1 \times S^2$ implies that $\chi(\HFp_{\leq
d(Y,\spinct)+2N+1}(Y, \spinct))$ is independent of $N$ for
sufficiently large $N$. We let $\chiTrunc(Y, \spinct)$ denote the
value of this Euler characteristic.

One can express $\chiTrunc$ in terms of Turaev torsion function
\cite{Turaev}
$$\tau_Y \colon \SpinC(Y) \longrightarrow \Z.$$
It is proved in \cite{HolDiskTwo} that  for any $\spinct$,
$$\chiTrunc(Y, \spinct)=-\tau(Y,\spinct).$$
For the precise statement and the sign issues for Turaev function
we refer to Proposition 10.14  of \cite{HolDiskTwo}.

In what follows, $\lambda$ denotes the \emph{Casson-Walker}
invariant normalized by $\lambda(\Sigma(2,3,5))=-1$, where
$\Sigma(2,3,5)$ is oriented as the boundary of the negative
definite $E_8$ plumbing.

\medskip
\noindent {\bf Surgery}\qua Here we set up our framework for
surgeries. We directly follow \cite{Theta}. Let $X$ be an oriented
three-manifold with a torus boundary and $H_1(X;\R)\cong \R$.  The
map $H_1(\partial X;\Z) \longrightarrow H_1(X;\Z)$ has
one-dimensional  kernel. Let $\ell'$ denote a generator for the
kernel, $d(X)>0$ denote its divisibility, and let $\ell$ be the
element $\ell'/d$. We call $\ell$ the {\em longitude}.

Fix a homology class $m\in H_1(\partial X)$ with $m\cdot \ell=1$.
For a pair of relatively prime integers $(p,q)$, the manifold
$Y_{p/q}$ is obtained from $X$ by attaching a $S^1\times \CDisk$
with $\partial \CDisk= p m + q \ell$, and let $Y=Y_{1/0}$. Note
that in general $Y_{p/q}$ depends on a choice of $m$, but
$Y_{0}=Y_{0/1}$ does not. Note also that $Y_{0}$ is a rational
homology $S^1\times S^2$, while all the other $Y_{p/q}$ are
rational homology spheres.

There is a short exact sequence
\begin{equation}
\label{eq:SESYnot}
\begin{CD}
0@>>> \Z @>>> \SpinC(Y_0) @>>> \SpinC(X) @>>> 0
\end{CD},
\end{equation}
by which we mean that the subgroup $\Z\subset H^2(Y_0;\Z)$
generated by the Poincar\'e dual to $m$ (viewed as a subset of
$Y_0$) acts freely on $\SpinC(Y_0)$, and its quotient is naturally
identified (under restriction to $X\subset Y_0$) with $\SpinC(X)$.

Thus, each $\SpinC$ structure $\spinc$ on $X$ has a natural level
$y=y(\spinc)\in\Zmod{d}$ defined as follows.  Let $\spincb$ be any
$\SpinC$ structure on $Y_0$ whose restriction is $\spinca$, and
consider its image in $$\tSpinC(Y_0)/\Z(\PD[m])\cong \Zmod{d},$$
where $\tSpinC(Y_0)$ is the group of $\SpinC$ structures modulo
the action of the torsion subgroup of $H^2(Y_0;\Z)$.

Furthermore, for any of the $Y_{p/q}$, the map $\SpinC(Y_{p/q})$
to $\SpinC(X)$ is surjective, and its fibers consist of orbits by
a cyclic group generated by the Poincar\'e dual to the knot which
is the core of the complement $Y_{p/q}-X$ (for $Y=Y_{1/0}$, this
fiber has order $d=d(X)$). For a fixed $\SpinC$ structure
$\spinca$ on $X$, let $\SpinC(Y_{p/q};\spinca)$ denote the set of
$\SpinC$ structures $\spincb\in\SpinC(Y_{p/q})$ whose restriction
to $X$ is $\spinca$.

\section{Surgery formula and its applications}
Our main theorem is the following surgery formula for the Euler
characteristic.
\begin{theorem}
\label{main} For integers $p,q,d,y$ with $p\neq0$, $p$ and $q$
relatively prime, $d>0$ and $0\leq y< d$, there is quantity
$\epsilon(p,q,d,y)\in\Q$ with the following property.  Let $X$ be
an oriented rational homology $S^1\times \CDisk$, with
divisibility $d(X)=d$, and choose $m$, $\ell$ as described in the
previous section.  Fixing any $\SpinC$ structure $\spinca$ over
$X$ with level $y(\spinca)=y$, we have the relation:
\begin{align*}
\sum_{\spincb\in\SpinC(Y_{p/q};\spinca)}\Eul(Y_{p/q},\spincb)=&
p\left(\sum_{\spincc\in\SpinC(Y;\spinca)}\Eul(Y,\spincc)\right)-q\left(\sum_{\spincd\in\SpinC(Y_0;\spinca)}\chiTrunc(Y_0,\spincd)\right)+\\
&\mbox{} + \epsilon(p,q,d,y).
\end{align*}

\end{theorem}
\begin{cor}
\label{cor:sum} For $X$ as above,
\begin{align*}
\left(\sum_{\spincb\in\SpinC(Y_{p/q})} \Eul(Y_{p/q},\spincb)
\right) = &p \left(\sum_{\spincc\in\SpinC(Y)}\Eul(Y,\spincc)
\right)+ q(\sum_{i=1}^\infty a_i i^2) \\ &\mbox{}+ |\Tors
H_1(X;\Z)| \epsilon(p,q,d),
\end{align*}
where $d=d(X)$, $a_i$ are the coefficients of the symmetrized
Alexander polynomial of $Y_0$, normalized so that
$$
A(1)=|\Tors H^2(Y_0;\Z)|,
$$
and
$$\epsilon(p,q,d)=\sum_{y=0}^{d-1}\frac{\epsilon(p,q,d,y)}{d}.$$
\end{cor}
\begin{proof}
This follows from the surgery formula and the fact that
$\chiTrunc(Y_0,\spinct) = -\tau(Y_0,\spinct)$.
\end{proof}

\begin{theorem}
\label{thm:CassonWalker} For any rational homology three-sphere
$M$ we have
$$\sum_{\spinct\in\SpinC(M)}\Eul(M,\spinct)=\big|H_1(M;\Z)\big|\lambda(M),$$ where $\lambda(M)$ is the
Casson-Walker invariant of $M$.
\end{theorem}
\begin{proof}
We already have a surgery formula for
$\sum_{\spinct\in\SpinC(Y)}\Eul(Y,\spinct)$. The scaled
Casson-Walker invariant
$$\lambda'(Y)=|H_1(Y;\Z)|\lambda(Y)$$
satisfies a similar formula with possibly different constants
$\epsilon'(p,q,d)$, see \cite{Theta}. In fact, we have
\begin{eqnarray*}\lambda'(Y_{p/q})= p \lambda'(Y) + q\left(\sum_{j\geq 1} a_j j^2\right) + |\Tors
H_1(X;\Z)| \left(\frac{q(d^2-1)}{24 d}-\frac{p d \cm
s(q,p)}{2}\right),
\end{eqnarray*}
i.e. $\epsilon'(p,q,d)= \left(\frac{q(d^2-1)}{24 d}-\frac{p d \cm
s(q,p)}{2}\right)$. Thus, it remains to show that
$$\epsilon(p,q,d)=\epsilon'(p,q,d).$$
For $d=1$ we can use a model calculation on $Y=S^3$ with the
surgery made on the unknot. Since $S^3_{p/q}= L(-p,q)$, by
\cite{Nemethi} (see also \cite{rasmussen}) we have
$$\sum_{\spinct\in\SpinC(L(p,q))}d(L(-p,q),\spinct)
=p\cm s(q,-p)=p\cm s(q,p).$$ Taking into account that
$\HFpred(L(-p,q)) \cong 0$ it follows that in this case
 $$\sum_{\spinct\in\SpinC(S^3_{p/q})}\Eul(S^3_{p/q},\spinct)=
-\frac{p\cm s(p,q)}{2}.$$ Plugging this into the surgery formula
we get $$\epsilon(p,q,1)= -\frac{p\cm s(p,q)}{2}$$ as needed.

To complete the proof, one shows that $\epsilon(p,q,d)$ is
determined by the surgery formula and the values of
$\epsilon(p,q,1)$. This is done by considering the Seifert
manifold $M(n,1,-n,1,q,-p)$. It can be  obtained from
$M(n,1,-n,1,0,1)$ by $(p,q,n)$ surgery. On the other hand, it is
possible to show that this manifold can be obtained by a sequence
of surgeries on knots with $d$'s less than $n$, see \cite{Theta}
for details.
\end{proof}
Now let us formulate the connection between the renormalized Euler
characteristic and Turaev torsion. For rational homology
three-sphere $M$ and a $\SpinC$ structure $\spinct$ on it define
$$\widehat{\tau}(M, \spinct) =-\tau(M,\spinct) +\lambda(M).$$

\begin{theorem}
For any rational homology three-sphere $M$ and a $\SpinC$
structure $\spinct$ on it we have
$$\Eul(M,\spinct)=\widehat{\tau}(M,\spinct)=SW(M,\spinct).$$
\end{theorem}
\begin{proof}
The proof follows using the framework of \cite{Nicola}. According
to it, there are several conditions on $\Eul$ and $\chiTrunc$ that
guarantee the sought equality. We list them as follows:
\begin{itemize}
\item The  surgery formula of Theorem \ref{main} is satisfied.
Note that we have a negative sign in front of the second term, but
it can be made positive by switching from $\chiTrunc$ to
$-\chiTrunc$. \item For any three-manifold $M$ with $b_1(M)=1$ and
a $\SpinC$ structure $\spinct$ on it
$$-\chiTrunc(M,\spinct) = \tau(M,\spinct).$$ \item  For any
rational homology sphere $M$,
$$\sum_{\spinct\in\SpinC(M)}\Eul(M,\spinct)=
\big|H_1(M;\Z)\big|\lambda(M).$$ \item For any integral homology
sphere $M$
$$\Eul(M,\spinct_0)=\widehat{\tau}(M,\spinct_0),$$
where $\spinct_0$ is the unique $\SpinC$ structure on $M$. \item
When $M$ is a Lens space
$$\Eul(M,\spinct)=\widehat{\tau}(M,\spinct),$$ for any $\SpinC$
structure $\spinct$ on $M$. \item If $M_1$ and $M_2$ satisfy
$\Eul=\widehat{\tau}$ then so does $M_1\#M_2$.\end{itemize} The
first three facts have already been mentioned, while the fourth
item is Theorem 5.1 of \cite{AbsGraded}, the fifth condition is
satisfied by \cite{Nemethi}. The last statement follows from
additivity of $d$, see Theorem 4.3 of \cite{AbsGraded} and from a
Kunneth type formula, see Corollary 6.3 of \cite{HolDiskTwo}. The
theorem follows.
\end{proof}

\section{Proof of the surgery formula}
Let  $\theta^c$ denote the three-dimensional $\SpinC$ homology
bordism group, defined as the set of equivalence classes of pairs
$(M,\spinct)$ where $M$ is a rational homology three-sphere, and
$\spinct$ is a $\SpinC$ structure over $M$, with the equivalence
given as follows. Ne say $(M_1,\spinct_1)\sim (M_2,\spinct_2)$ if
there is a (connected, oriented, smooth) cobordism $N$ from $M_1$
to $M_2$ with $H_i(N,\Q)=0$ for $i=1$ and $2$, which can be
endowed with a $\SpinC$ structure $\spinc$ whose restrictions to
$M_1$ and $M_2$ are $\spinct_1$ and $\spinct_2$ respectively. The
connected sum operation makes this set an Abelian group (whose
unit is $S^3$ with its unique $\SpinC$ structure).
 The invariant $d(M,\spinct)$ gives a group homomorphism
$$d \colon \theta^c \longrightarrow \Q.$$

It is proved in \cite{AbsGraded} that $d$ is a lift of the
classical homomorphism
$$\rho \colon \theta^c \longrightarrow \Q/2\Z$$
(see \cite{atiyah}) defined as follows. Let $N$ be any
four-manifold equipped with a $\SpinC$ structure $\spinc$ with
$\partial N \cong M$ and $\spinc|\partial N \cong \spinct$ then
$$\rho(M,\spinct)\equiv \frac{c_1(\spinc)^2-\sgn(N)}{4}
\pmod{2\Z}$$ where $\sgn(N)$ denotes the signature of the
intersection form of $N$.

Going back to our surgery notation, let $W$ be  the standard
cobordism between $Y$ and $Y_{p/q}$ obtained by 2-handle
additions. Let $\rho'(Y,\spinct) \equiv \rho(Y, \spinct)
\pmod{2\Z}$ such that $\rho'(Y,\spinct) \in [0,2)$. For the
manifold $Y_{p/q}$ and a $\SpinC$ structure $\spinct$ on it
consider any $\spinc$ on $W$ with $\spinc|Y_{p/q} = \spinct$. We
define $\rho'(Y_{p/q}, \spinct)=\rho'(Y, \spinc|Y)$.

For any constant $k$ define
$$\HFp_{\preceq k}(Y_{p/q}, [\spinca])= \bigoplus_{\spinct\in\SpinC(Y_{p/q};\spinca)} \bigoplus_{\{d\in\Q\big| d\leq k+\rho'(Y_{p/q},\spinct)\}}\HFp_d(Y_{p/q},\spinct).$$

$Y_0$ is not a rational homology sphere, if $\spinct$ is torsion
$\SpinC$ structure on $Y_0$ one can still define
$\rho'(Y_0,\spinct)$ similarly to above. It is useful to note that
equivalence $$d(Y_0,\spinct) \equiv 1+\frac{c_1(\spinc)^2 +
\sgn(W)}{4} +\rho'(Y_0,\spinct) \pmod{2\Z}$$ holds for any $\spinc
\in \SpinC(W)$ satisfying $\spinc|Y_0=\spinct$, this follows from
the grading shift formula for maps induced by cobordisms. One
should look at both absolute $\Q$ and $\Z/2\Z$ grading shifts.

Let $\spincset$  be the subset of torsion $\SpinC$ structures of
$\SpinC(Y_0)$. Now set $$\HFp_{\preceq k}(Y_0,
[\spinca])=\bigoplus_{\spinct\in \SpinC(Y_0;\spinca) \setminus
\spincset}\HFp(Y,\spinct) \oplus \bigoplus_{\spinct \in
\SpinC(Y_0;\spinca) \cap \spincset} \bigoplus_{\{d\in\Q\big| d\leq
k+\rho'(Y_0,\spinct \}} \HFp_d(Y_0,\spinct) .$$

\begin{lemma}
\label{EulerRed} For integers $p,q,d,y$ with $p$ and $q$
relatively prime, $d>0$ and $0\leq y< d$, there is quantity
$k(p,q,d,y)$ with the following property.  Let everything be as in
Theorem~\ref{main}, then
\begin{eqnarray*}
\lefteqn{\chi(\HFp_{\preceq 2N}(Y_{p/q},[\spinca]) - N \cm
|\SpinC(Y_{p/q};\spinca)|=}
\\
& &=\sum_{\spincb\in\SpinC(Y_{p/q};\spinca)}\Eul(Y_{p/q},\spincb)+
p \sum_{\spincc\in\SpinC(Y;\spinca)} \frac{\rho'(Y,\spincc)}{2}+
k(p,q,d,y).
\end{eqnarray*}
\end{lemma}
\begin{proof}
(cf. lemma 4.8 in \cite{unk}.) For sufficiently large $N$,
$\HFpRed(Y_{p/q},[\spinca])$ is contained in $\HFp_{\preceq
2N}(Y_{p/q},[\spinca])$. Over $\Z$, we have a splitting
$$
\HFp_{\preceq 2N}(Y_{p/q},[\spinca]) \cong
\HFpRed(Y_{p/q},[\spinca])\oplus (\Image\pi \cap \HFp_{\preceq
2N}(Y_{p/q},[\spinca])).
$$
But it follows readily from the structure of $\HFinf(Y_{p/q})$
(c.f. Equation~\eqref{eq:HFinf}) that
\begin{align*}
\chi(\Image\pi\cap & \HFp_{\preceq 2N}(Y_{p/q},[\spinca])) =\\
&=\sum_{\spincb\in\SpinC(Y_{p/q};\spinca)}
\#\{{[d(Y_{p/q},\spincb),2N+\rho'(Y_{p/q},\spincb)]}\cap (d(Y_{p/q},\spincb)+2\Z)\subset \Q\} \\
&= \sum_{\{\spincb\in\SpinC(Y_{p/q};\spinca)\}} \left(N + 1 -\left
\lceil \frac{d(Y_{p/q},\spincb)-\rho'(Y_{p/q},\spincb)}{2}\right
\rceil\right),
\end{align*}
where here $\lceil x\rceil$ denotes the smallest integer greater
than or equal to $x$.  Thus we get that
\begin{align*}
\chi(\HFp_{\preceq
2N}&(Y_{p/q},[\spinca]))-N\cm|\SpinC(Y_{p/q};\spinca)|=
\\&=\sum_{\spincb\in\SpinC(Y_{p/q};\spinca)} \left(\HFpRed(Y_{p/q},\spincb)-\left \lceil\frac{d(Y_{p/q},\spincb)-\rho'(Y_{p/q},\spincb)}{2}\right \rceil+1\right)
\\&=\sum_{\spincb\in\SpinC(Y_{p/q};\spinca)}\left(\HFpRed(Y_{p/q},\spincb)-\frac{d(Y_{p/q},\spincb)}{2}\right)
+\\&+\sum_{\spincb\in\SpinC(Y_{p/q};\spinca)}\left(\frac{d(Y_{p/q},\spincb)}{2}-\left
\lceil\frac{d(Y_{p/q},\spincb)-\rho'(Y_{p/q},\spincb)}{2}\right
\rceil+1\right).
\end{align*}

To complete the proof we have to show that the difference

$$k=\sum_{\spincb\in\SpinC(Y_{p/q};\spinca)}
\left(\frac{d(Y_{p/q},\spincb)}{2}-\left
\lceil\frac{d(Y_{p/q},\spincb)- \rho'(Y_{p/q},\spincb)}{2}\right
\rceil+1\right) - p\cm
\sum_{\spincc\in\SpinC(Y;\spinca)}\frac{\rho'(Y,\spincc)}{2}$$
depends only on $p,q,d,y$. Clearly
$$k=\sum_{\spincb\in\SpinC(Y_{p/q};\spinca)}
\left(\frac{d(Y_{p/q},\spincb)-
\rho'(Y_{p/q},\spincb)}{2}-\left\lceil\frac{d(Y_{p/q},\spincb)-
\rho'(Y_{p/q},\spincb)}{2}\right\rceil+1\right).$$ This in turn
depends only on $d(Y_{p/q},\spincb)- \rho'(Y_{p/q},\spincb)
\pmod{2\Z}=\rho(Y_{p/q},\spincb)- \rho'(Y_{p/q},\spincb)
\pmod{2\Z}$ which is completely determined by the collection of
all $c_1(\spinc)^2 \pmod{8\Z}$ with $\spinc \in \SpinC(W)$
satisfying $\spinc|Y \in \SpinC(Y;\spinca)$. This follows from the
definitions and the fact that $\rho$ is a homomorphism. Hence, the
proof is concluded by the following lemma.

\begin{lemma}
\label{collection} Let $W$ be the standard cobordism between $Y$
and $Y_{p/q}$. The collection with repetitions of all
$c_1(\spinc)^2$ satisfying $\spinc \in \SpinC(W)$ and $\spinc|Y
\in \SpinC(Y;\spinca)$ is completely determined by the values of
$p, q, d$ and $y$.\qed
\end{lemma}
\end{proof}
\begin{lemma}
For integers $d,y$ with $d>0$ and $0\leq y< d$, there is quantity
$r(d,y)$ with the following property.  Let everything be as in
Theorem~\ref{main}, then
$$\chi(\HFp_{\preceq 2N}(Y_{0},[\spinca])= \sum_{\spincb\in\SpinC(Y_{0};\spinca)}\chiTrunc(Y_{0},\spincb)+
r(d,y).$$
\end{lemma}
\begin{proof}
The idea of the proof is the same with the previous one. We do not
have any terms involving $N$ because of the different structure of
$\HFinf$ for manifolds with $b_1=1$.
\end{proof}

\begin{lemma}
\label{Euler} For integers $p,q,d,y$ with $p\neq 0$, $p$ and $q$
relatively prime, $d>0$ and $0\leq y< d$, there is quantity
$c(p,q,d,y)$ with the following property.  Let everything be as in
Theorem~\ref{main}, then
\begin{eqnarray}
\chi (\HFp_{\preceq 2N}(Y_{p/q},[\spinca])) = p\cm
\chi(\HFp_{\preceq 2N} (Y,[\spinca])) - q \cm \chi(\HFp_{\preceq
2N}(Y_0,[\spinca])) + c(p,q,d,y),
\end{eqnarray}
provided that $N$ is sufficiently large.
\end{lemma}
\begin{proof}
The proof is a generalization of the argument of lemma 4.9 in
\cite{unk}. Let us use induction on $p+q$. The base of induction
is the case when $p+q=1,2$, which reduces to $(p,q)=(1,0)$ or
$(1,1)$. The lemma clearly holds for the first combination; we
will discuss the second case in the end of the proof.

For a pair $(p,q)$ of relatively prime, non-negative integers with
$p+q>2$, one can select two pairs of non-negative, relatively
prime integers $(p_0,q_0)$ and $(p_2,q_2)$, with $p_0, p_2\neq 0$
satisfying
\begin{eqnarray}
p_0 \cm q - p \cm q_0 &=& -1 \\
(p,q)&=& (p_0,q_0)+(p_2,q_2)
\end{eqnarray}
Consider the manifolds $\Yzero$, $\Yone$ and $\Ytwo$. There are
standard 2-handle cobordisms between these manifolds. Let $W_0$
denote the cobordism between $\Yzero$ and $\Yone$, $W_1$ the
cobordism between $\Yone$ and $\Ytwo$, $W_2$ between $\Ytwo$ and
$\Yzero$. We can write down the following long exact sequence
$$\begin{CD} ...@>>> \HFp(\Yzero,[\spinca])
@>{f_0}>>\HFp(\Yone,[\spinca]) @>{f_1}>>\HFp(\Ytwo,[\spinca])
@>{f_2}>> ...,
\end{CD}$$
where the maps are induced by the corresponding cobordisms. Note
that $W_0$ and $W_1$ are both negative definite, but $W_2$ is not.

By inductive hypothesis the lemma holds for $(p_0,q_0)$ and
$(p_2,q_2)$. When $N$ is sufficiently large, the image of the
restriction $g_0$ of $f_0$ to $\HFp_{\preceq
2N}(\Yzero,[\spinca])$ is contained in $\HFp_{\preceq
2N+\frac{1}{4}}(\Yone,[\spinca])$, the restriction $g_1$ of $f_1$
to $\HFp_{\preceq 2N+\frac{1}{4}}(\Yone,[\spinca])$ is contained
in $\HFp_{\preceq 2N+\OneHalf}(\Ytwo,[\spinca])$, and finally, the
restriction $g_2$ of $f_2$ to $\HFp_{\preceq
2N+\OneHalf}(\Ytwo,[\spinca])$ is contained in $\HFp_{\preceq
2N}(\Yzero,[\spinca])$. This follows at once from the definition
of $\rho'$ which appears in the expression for $\HFp_{\preceq k}$,
and the grading shift formula: we have that $\chi(W_i)=1$ and
$\sigma(W_i)=-1$ for $i=0,1$; while the cobordism $W_2$ induces
the trivial map on $\HFinf$ since $b_2^+(W_2)=1$.

Choosing $N$ as above, consider the diagram

\begin{scriptsize}
$$
\begin{CD}
&& 0 &&  0 && 0  \\
&& @VVV @VVV @VVV \\
... @>{g_2}>> \HFp_{\preceq 2N}(\Yzero,[\spinca]) @>{g_0}>> \HFp_{\preceq 2N+\frac{1}{4}}(\Yone,[\spinca]) @>{g_1}>> \HFp_{\preceq 2N+\OneHalf}(\Ytwo,[\spinca]) @>{g_2}>> ... \\
&& @VVV @VVV @VVV \\
... @>{f_2}>> \HFp(\Yzero,[\spinca]) @>{f_0}>>\HFp(\Yone,[\spinca]) @>{f_1}>> \HFp(\Ytwo,[\spinca]) @>{f_2}>> ... \\
&& @VVV @VVV @VVV \\
... @>{h_2}>> \HFp_{\succ 2N}(\Yzero,[\spinca])@>{h_0}>>\HFp_{\succ 2N+\frac{1}{4}}(\Yone,[\spinca]) @>{h_1}>>\HFp_{\succ 2N+\OneHalf}(\Ytwo,[\spinca])  @>{h_2}>> ..., \\
&& @VVV @VVV @VVV \\
&& 0 &&  0 && 0  \\
\end{CD}
$$
\end{scriptsize}
where the columns are exact. Note that the first and the third
rows are not necessarily exact, while the middle one is exact. Let
us think of these three rows as chain complexes. We denote these
three rows by ${\mathcal R}_1$, ${\mathcal R}_2$, and ${\mathcal
R}_3$. Since ${\mathcal R}_2$ is exact, it follows that
$H_*({\mathcal R}_1)\cong H_*({\mathcal R}_3)$.

Now let us show that $H_*({\mathcal R}_3)$ is determined by
$p,q,d$ and $y$ for $N$ sufficiently large. This is established
using the structure of maps on $\HFinf$, lemma \ref{collection}
and the diagram

\begin{scriptsize}
$$
\begin{CD}
... @>{h_2^{\infty}}>> \HFinf_{\succ 2N}(\Yzero,[\spinca])@>{h_0^\infty}>>\HFinf_{\succ 2N+\frac{1}{4}}(\Yone,[\spinca]) @>{h_1^\infty}>>\HFinf_{\succ 2N+\OneHalf}(\Ytwo,[\spinca])  @>{h_2^{\infty}}>> ... \\
&& @V{\cong}VV @V{\cong}VV @V{\cong}VV \\
... @>{h_2}>> \HFp_{\succ 2N+}(\Yzero,[\spinca])@>{h_0}>>\HFp_{\succ 2N+\frac{1}{4}}(\Yone,[\spinca]) @>{h_1}>>\HFp_{\succ 2N+\OneHalf}(\Ytwo,[\spinca])  @>{h_2}>> ..., \\
\end{CD}
$$
\end{scriptsize}

where here $h_0$ is the sum over all $\spinc\in\SpinC(W_0)$ of the
projections of the induced maps on $\HFinf$; e.g. letting
$$\Pi_{\succ 2N+\OneHalf}\colon \HFinf(\Yone,[\spinca]) \longrightarrow \HFinf_{\succ 2N+\OneHalf}(\Yone,[\spinca])$$
denote the projection, we let $h_0^{\infty}$ be the restriction to
$\HF_{\succ 2N}(\Yzero,[\spinca])$ of
$$\sum_{\spinc\in\SpinC(W_0)} \Pi_{\succ 2N+\OneHalf}\circ F^{\infty}_{W_0,\spinc}.$$The maps $h_i^{\infty}$ are defined similarly. Note that
$h_2^{\infty}=0$, since the map induced by $W_2$ has
$b_2^+(W_2)=1$.

So far we have established that for all sufficiently large $N$,
$$\chi (H_*({\mathcal R}_1)) = \chi (\HFp_{\preceq 2N}(\Yzero,[\spinca])) - \chi(\HFp_{\preceq 2N+\OneQuarter}(\Yone,[\spinca]))
+\chi(\HFp_{\preceq 2N+\OneHalf}(\Ytwo,[\spinca]))$$ is completely
determined by $p, q, d, y$. It is also clear that for sufficiently
large $N$,
$$
\chi(\HFp_{\preceq 2N+\OneQuarter}(\Yone,[\spinca]))=\chi
(\HFp_{\preceq 2N}(\Yone,[\spinca])) + c_3$$ and
$$\chi(\HFp_{\preceq 2N+\OneHalf}(\Ytwo,[\spinca]))=\chi (\HFp_{\preceq 2N}(\Ytwo,[\spinca])) + c_4,$$
with constants $c_3$ and $c_4$ again depending only on $p_0,q_0,
d, y$ and $p_2,q_2, d, y$ respectively, lemma \ref{collection}.
Combining all the constants, we establish the inductive step in
the case where $p_0$ is non-zero.

When  $(p,q)=(1,1)$, the above argument works with slight
modification. In this case, we consider the manifolds $Y, Y_0,
Y_1$. The dimension shifts work differently:
$\sigma(W_0)=\sigma(W_1)=0$ and hence, we compare $\HFp_{\preceq
2N}(Y,[\spinca])$, $\HFp_{\preceq 2N+\OneHalf}(Y_0,[\spinca])$,
and $\HFp_{\preceq 2N+1}(Y_1,[\spinca])$. To see that the map
$f_2$ induced by $W_2$ carries $\HFp_{\preceq
2N+1}(Y_1,[\spinca])$ into  $\HFp_{\preceq 2N}(Y,[\spinca])$ for
sufficiently large $N$, remember that the kernel of the map $f_0$
induced by $W_0$ is finitely generated. Some parities change under
these maps, so the Euler characteristic is given as follows
$$\chi (H_*({\mathcal R}_1)) = \chi (\HFp_{\preceq 2N}(Y,[\spinca]))
- \chi(\HFp_{\preceq 2N+\OneQuarter}(Y_0,[\spinca]))
-\chi(\HFp_{\preceq 2N+\OneHalf}(Y_1,[\spinca])),$$ compare with
Proposition 5.3 in \cite{AbsGraded}.

\end{proof}
\noindent{\bf{Proof of Theorem 5.3.1}}\qua When $p$ and $q$ are
non-negative, this is a combination of Lemmas~\ref{EulerRed} and
\ref{Euler}. The remaining case can be proved by running the
induction from Lemma~\ref{Euler} to show that it still holds in
the case where $p>0$ and $q\leq 0$. \qed

\end{document}